\documentstyle[11pt]{article}
\newfont{\bbold}{msbm10 scaled \magstep1}
\newfont{\bbolds}{msbm7 scaled \magstep1}

\newcommand{\N}{\mbox{\bbold N}}
\newcommand{\Z}{\mbox{\bbold Z}}
\newcommand{\Q}{\mbox{\bbold Q}}

\newcommand{\C}{\mbox{\bbold C}}

\newcommand {\ord}{\mathop{\rm ord}}

\newcommand {\re}{\mathop{\rm Re}}

\newcommand {\inop}{\in \C [z,E]}

\newcommand {\e }{\varepsilon }
\newcommand {\f }{\varphi }

\newcommand {\k }{\kappa  }

\newcommand {\qed}{\hfill$\Box$\par\medskip}

\newtheorem{defi}{Definition}
\newtheorem{rem}{Remark}
\newtheorem{Theorem}{Theorem}
\newtheorem{Proposition}{Proposition}

\newtheorem{Example}{Example}

\title {
Apparent Singularities
of Linear
Difference  Equations
with Polynomial Coefficients
}
\author
{
\large\sl S. A.~Abramov\thanks{Supported by
French-Russian Lyapunov Institute  grant 98-03
and by RFBR grant 04-01-00757. } \\
\normalsize Russian Academy of Sciences \\
\normalsize  Dorodnicyn Computing Centre\\
\normalsize  Vavilova 40, 119991, Moscow GSP-1\\
\normalsize Russia\\
\normalsize {\tt sabramov@ccas.ru}
\and \large\sl M. A.~ Barkatou \\
\normalsize  Universit\'e de Limoges,   LACO\\
\normalsize 123, Av.  A. Thomas\\
\normalsize 87060 Limoges cedex  \\
\normalsize France\\
\normalsize {\tt moulay.barkatou@unilim.fr}
\and \large\sl M. van Hoeij\thanks{Supported by
NSF grant 0098034.} \\
\normalsize Florida State University\\
\normalsize Department of  mathematics \\
\normalsize Tallahassee, FL 32306-3027 \\
\normalsize USA\\
\normalsize {\tt hoeij@math.fsu.edu}
}

\begin{document}
\date{}

\maketitle
\begin{abstract}

Let $L$ be a linear difference operator with polynomial coefficients.
We consider singularities of $L$ that correspond to roots of the
trailing (resp. leading) coefficient of $L$.
We prove that
one can effectively construct
a left multiple with polynomial coefficients ${\tilde L}$ of $L$
such that
every singularity of $\tilde{L}$ is
a singularity of $L$ that is not apparent.
As a consequence,
if all singularities of $L$ are apparent, then $L$ has
a left multiple whose trailing  and leading coefficients equal 1.
\end{abstract}
\section{Introduction}
\label{intro}
Investigation of singular points (also called singularities) of a
linear differential or difference operator $L$ gives an
opportunity to study singularities of
solutions of the equation $L(y)=0$ without solving this
equation.

In the differential case, take $L = \sum_{i=0}^d a_i(z) \partial^i \in \C[z,\partial]$
where $d$ is the order of $L$ (so $a_d(z) \neq 0$) and $\partial = d/dz$.
Let $S(L) \subset \C$ be the set of roots of $a_d(z)$.
The finite singularities of $L$ are the elements of $S(L)$. The
singularity at $\infty$ will not be considered in this paper.
If $p \in S(L)$ and if
there exist $d$ linearly independent analytic solutions at $z=p$
then $p$ is called an {\em apparent} singularity.
Suppose $L$ has apparent singularities. The question is if
it is possible to
construct another operator ${\tilde L} \in \C[z,\partial]$ of higher order such that any
solution of  $L(y)=0$ is a solution of ${\tilde L}(y)=0$,
and $S({\tilde L}) = \{ p \in S(L) \ | \ p {\rm \ not \ apparent} \}$.
In the differential case the answer is affirmative, see \cite{Tsai}.
In this paper we give
the affirmative answer to the corresponding question for the difference case.

In the remainder of this paper (except the appendix) only the
difference case will be considered.
The shift operator $E$ acts on functions of  the complex variable $z$
as $Ey(z)=y(z+1)$. We consider non-commutative operator rings
$\C[z,E]$ and $\C(z)[E]$ (the rings of linear
difference operators with polynomial and, resp., rational
function coefficients over $\C$).
Let
\begin{equation}
\label{eq10}
L=a_d(z)E^d+\cdots +a_1(z)E+ a_0(z)
\in \C[z,E].
\end{equation}
Assume that the leading
coefficient $a_d(z)$ and the trailing coefficient $a_0(z)$ are both
non-zero,
and that $a_0(z),\dots ,a_d(z)$ do not have a non-constant common
factor.
Set $\ord L=d$.
\begin{defi}
\label{tls}
A root $p$ of $a_0(z)$
is called a
$t$-{\em singularity}
(a trailing singularity).
A root $p$ of $a_d(z-d)$ is called an
$l$-{\em singularity}
(a leading singularity).
\end{defi}
\begin{defi}
\label{halfh}
A {\em right-holomorphic} (resp. {\em left-holomorphic}) function
is a meromorphic function on $\C$
that is holomorphic on some
right (resp. left) half plane. In other words, holomorphic
when $\re z$ (resp. $-\re z$) is sufficiently large.
A {\em half-holomorphic} function is a function that is
right- or left-holomorphic.
\end{defi}

\begin{defi}
\label{defapp}
A root $p$ of $a_0(z)$ (resp. of $a_d(z-d)$)
is called an {\em apparent}
$t$-(resp. $l$)-singularity if
no right-(resp. left)-holomorphic solution has a pole at $p$.
An operator ${\tilde L}$ is a
$t$-(resp. $l$)-desingularization of $L$ if
every meromorphic solution of $L$ is a solution of $\tilde{L}$,
and every $t$-(resp. $l$)-singularity of $\tilde{L}$ is
a $t$-(resp. $l$)-singularity of $L$ that is not apparent.
\end{defi}

We show that both $t$- and $l$-desingularizations exist.
We give algorithms {\tt t-desing} and {\tt l-desing} for
constructing a $t$-(resp. $l$)-desingularization and algorithm {\tt
desingboth} for constructing a desingularization related to both
trailing and leading coefficients.

The above definition of a desingularization is not the same as in
\cite{AvH99}
(see the summary of \cite{AvH99} given in Section~\ref{42}).

Our approach is based on some specific properties
of apparent
singularities that are proved in this paper
(Propositions \ref{propn},\ref{propn'}).

Besides of a theoretical interest, it is useful to have
a desingularization ${\tilde L}$
of $L$ for solving the continuation problem.
Equation $L(y)=0$ can be used as a tool to
define a sequence or a function. If we know the value $y(z)$ at
every point $z$ of a  given strip $ \lambda \leq \re  z < \lambda
+\delta$, where $\delta$ is larger or equal to the order of $L$
then we can  find the value of $y(z)$ in the strip $
\lambda -1 \leq \re  z < \lambda$, and then in the strip $ \lambda -2
\leq \re  z < \lambda-1$, and so on. We can keep continuing $y(z)$ to
the left in this way except when we encounter $t$-singularities.
Similarly, we can
continue $y(z)$ to the right except at $l$-singularities.
Thus, the singularities of $L$ may present obstacles to continuing
solutions of $L(y)=0$.
If the singularities are apparent, one can overcome those obstacles
using
$\tilde{L}$ instead of $L$.
Another use of desingularization is the following:
In the process of continuing sequences to the left (resp. to the right)
one must always divide by the trailing (resp. leading) coefficient.
If all $t$-(resp. $l$)-singularities are apparent, then
one can avoid such divisions by computing a desingularization, see
Sect.~\ref{appds} for an application.
However, there is a price to pay, namely that the order increases.

In the appendix we give algorithm $\partial$-{\tt desing}
for desingularization in the differential case.
More general results  for the differential case can 
be found in \cite{Tsai}, where $\C[z,\partial] \bigcap
\C(z)[\partial]L$ is computed. 
We include this appendix for completeness and because the proof is
short. \\

\noindent {\bf Acknowledgment.} The authors thank the anonymous
referees for numerous valuable and helpful comments.

\section{The sets
$C_{q,\sigma}(L)$ and $R_{q,\sigma}(L)$}
\subsection{
The set of singularities}
We consider  a linear difference operator (\ref{eq10}).
The set of trailing resp. leading
singularities is
$S_t(L) = \{p \in \C \,|\,a_0(p)=0\}$
resp. $S_l(L) = \{p \in \C \,|\,a_d(p-d)=0\}$.
The set of singularities is $S(L) = S_t(L) \bigcup S_l(L)$.

A point $p \in \C$ is
said to be {\em congruent} to a $t$-(resp. $l$)-singularity of $L$ if
$a_0(p+\nu)=0$ (resp. $a_d(p-d-\nu)=0$)
for some $\nu\in \N$.
Let $L$ be of the form
(\ref{eq10}) and
$\sigma_1,\dots, \sigma _m$ be all singularities of $L$.
Set
         \begin{equation}
         \label{eq10''}
         \iota (L)=\min \{\infty, \re \sigma_1,\dots ,\re \sigma _m\},\;
         \kappa (L)=\max \{-\infty, \re \sigma_1,\dots ,\re \sigma _m\}.
         \end{equation}
So $L$ has no singularity in the half-planes $\re z<\iota (L)$
and $\re z>\kappa (L)$.

Note that
any solution $F(z)$ of
$L(y)=0$ which is defined and holomorphic on a half-plane
$\re z>\k(L)$ (or, resp., $\re z<\iota (L)$),
can be continued to a meromorphic solution defined on $\C$
whose poles are congruent to $t$-(resp. $l$)-singularities.

Starting from this point until the end of Section \ref{sectionds}
we will consider only the $t$-singularities (the
$l$-singularities can be handled similarly).

\subsection{Systems of linear relations}
\label{ep}
Let $\Omega$ be a non-empty open subset of $\C$ which is stable under $E$,
i.e.,  $z \in \Omega$ implies  $z +1 \in \Omega$.
For example let $\Omega $ be the right half-plane $\re z > \k(L)$.
Let $\f  : \Omega \rightarrow \C$ be an arbitrary holomorphic function.
We associate to $\f$ a new function $\hat {\f}$ whose values are
formal Taylor power series in $\e$:  $$ \hat {\f}\,:\, \Omega
\longrightarrow \C[[\e ]], \ \ \
\hat{\f}(z) =
\sum_{\nu=0}^{\infty}{\frac{\f^{(\nu)}(z)}{\nu !} \e^\nu}.$$ Here
$\e$ is a new variable, rather than a ``small number''. Of course,
when $\e \in \C$ with $\vert \e \vert$ small enough the formal
series $\hat{\f}(z)$ converges and its sum is equal to $\f(z+\e)$.

If $\f$ is a polynomial then we can identify $\hat{\f}(z) = \f(z+\e)$.
The operator $L\in \C [z,E]$ has polynomial coefficients. We associate
to $L$ the operator
$$
\hat {L}=\hat{a}_d(z)E^d+\cdots +\hat{a}_0(z) = a_d(z+\e)E^d + \cdots + a_0(z+\e)
$$
which acts on functions $\Phi(z)$ whose values
are formal power series in $\e$.

The operator $\hat{L}$ acts also on sequences with values in
the field $\C((\e))$
of formal Laurent series.
If a finite sequence $f_q, f_{q+1}, \dots,f_{q+d-1} \in
\C((\e))$ is given for some $q\in \C$, then, by using the recurrence
given by $\hat{L}$, one can compute series \begin{equation} \label{4}
f_{q-1},f_{q-2},\dots \end{equation} An advantage of $\hat{L}$ in
comparison with $L$ is that neither the leading nor the trailing
coefficient of $\hat{L}$ vanishes when $z$ is any complex number.
However, a series $f_{q-m}$ can turn out to be
formal Laurent series for some positive integer $m$ even
when $f_q,\dots,f_{q+d-1}$ are formal Taylor series;
if $a_0(z_0)=0$ then the series $\hat{a}_0(z_0)$ is not invertible
in the ring $\C[[\e]]$.
The problem of defining complex values of
solutions at singularities of $L$ is transformed into the
problem of obtaining formal
Taylor
series (truncated power
series in computations).

Let $\sigma$ be some $t$-singularity of $L$. Now choose $q \in \sigma + \N$
for which $\re q > \k(L)$.
Let $\Phi (z)$ be a function
whose values are formal power series in $\e$,  and suppose that
$\hat{L} (\Phi)=0$.
If the values of $\Phi(z)$ at points $q,q+1,\dots, q+d-1$
are formal Taylor series
\begin{equation} \label{5}
\begin{array}{ccccccccc}
\Phi(q) &=& F_{q,0}&+&F_{q,1}\e&+&F_{q,2}\e ^2&+&\cdots\\
\vdots&&\vdots&&\vdots&&\vdots&&\\
\Phi(q+d-1) &=& F_{q+d-1,0}&+&F_{q+d-1,1}\e&+&F_{q+d-1,2}\e ^2&+&\cdots
\end{array}
\end{equation}
with $F_{ij} \in \C$
then using the equality $\hat {L}(\Phi)=0$
we can compute the formal
Laurent series $\Phi(\sigma)$, and each coefficient  of this
series will be a linear form in a finite set of $F_{ij}$'s.
This series $\Phi(\sigma)$ can contain negative exponents of $\e$.
We can find conditions on the coefficients $F_{ij}$'s in (\ref{5})
that guarantee that $\Phi (\sigma)$ is a Taylor series.
Indeed, if we use
the generic power series (\ref{5}) for this computation and if we get
some terms with negative exponents of $\e$, then after equating their
coefficients
to zero, we get a system of linear relations. This system forms
a necessary and sufficient condition that initial conditions (\ref{5})
lead to a
formal Taylor series $\Phi(\sigma)$ when $\Phi$ satisfies
$\hat{L}(\Phi)=0$.
This gives a finite system of
linear relations, denoted as $C_{q,\sigma}(L)$.  Constructing
$C_{q,\sigma}(L)$ can be carried out algorithmically, see also
\cite{AvH03} and \cite{vh}.  The system $C_{q,\sigma}(L)$ is a
necessary and sufficient condition that a function $\Phi(z)$, having
values (\ref{5}) at $q,\ldots,q+d-1$ and satisfying $\hat{L}(\Phi)=0$
has formal Taylor series values at $z=\sigma$.

\begin{Example}
\label{ex1}
For the operator
$$(z-1)zE^2-(3z+7)(z-3)E+(z+2)(z+1)$$
we can take generic formal Taylor series
$$
\begin{array}{ccccccc}
\Phi(4) = F_{40}&+&F_{41}\e&+&F_{42}\e ^2&+&\cdots\\
\Phi(5) = F_{50}&+&F_{51}\e&+&F_{52}\e ^2&+&\cdots
\end{array}
$$
Then the system $C_{4,-1}(L)$ is $\{20F_{40}-39F_{50}=0\}$.
This means that if $\Phi(4)$ and $\Phi(5)$ are as above, and if $\Phi$
is a solution of $\hat{L}$, then $\Phi(-1)$ is a Taylor series
if and only if $20F_{40}-39F_{50}=0$.
\end{Example}

\begin{Proposition}
\label{prop3}
With notations as above, if $C_{q,\sigma}(L)\neq \emptyset$
then $C_{q,\sigma}(L)$ contains
a relation $f(F_{q,0},F_{q+1,0},\dots, F_{q+d-1,0})=0$, where $f$ is a
non-zero linear form in variables  $F_{q,0},F_{q+1,0},\dots,
F_{q+d-1,0}$.
\end{Proposition}
\noindent {\bf Proof:}
Write $\Phi(\sigma) = \sum_{i=-N}^{\infty} F_{\sigma,i} \e^i$
with $F_{\sigma,-N}$ not zero. Now $N>0$ since $C_{q,\sigma}(L)\neq
\emptyset$. Furthermore, $N$ is the highest pole order of
any solution of $\hat{L}$ that has
formal Taylor
series as initial values at the points $q,q+1,\dots, q+d-1$.

Now $F_{\sigma,-N}$ can be written
as a linear combination of the $F_{q+i,j}$ with $0 \leq i < d$
and $0 \leq j$. Suppose that some $F_{q+i,j}$ with $j>0$ appears
with non-zero
coefficient.
Now assign the value 1 for this $F_{q+i,j}$
and the value 0 for all $F_{q+i',j'}$ with $(i',j') \neq (i,j)$.
Then we have obtained a function $\Phi$ with initial
values in $\{0, \e^j\}$ at $q,q+1,\dots, q+d-1$
and a pole of order $N$
at $\sigma$. Dividing by $\e^j$, we get initial values
in $\{0,1\} \subset \C[[\e]]$ and a pole of order $N+j$ at $\sigma$
which is a contradiction since $N+j>N$.
\qed

\subsection{The set $R_{q,\sigma}(L)$}
Let $\sigma$ be some $t$-singularity of $L$ and $q \in \sigma + \Z$. 
Let $\Phi (z)$ be a function
whose values are
rational functions of $\e$ (we do not expand them into power series
yet),  and suppose that $\hat{L}
(\Phi)=0$.  If the values of $\Phi(z)$ at points $q,q+1,\dots, q+d-1$
are known rational functions of $\e$
\begin{equation}
\label{phiphi}
\Phi (q), \Phi (q+1),\dots ,\Phi (q+d-1),
\end{equation}

then we can compute step-by-step the rational functions
\begin{equation}
\label{phiphi'}
\ldots, \Phi (q-2), \Phi (q-1),\dots ,\Phi (q+d), \Phi(q+d+1), \ldots
\end{equation}
in particular $\Phi(\sigma)$ using $\hat{L}$.
Consider the following $d$-tuples
\begin{equation}
\label{dtup}
(1,0,0,\dots, 0),(0,1,0,\dots, 0),
\dots
(0,0,0,\dots, 1)
\end{equation}
as  $d$ sets of initial values (\ref{phiphi}). For each set of
initial values one obtains a sequence of the form (\ref{phiphi'})
so we get rational functions 
$\Phi _i (\sigma), \;\;i=1,\dots ,d.$
Let
$$R_{q,\sigma}(L)= \{
\Phi _1 (\sigma), \dots, \Phi _d (\sigma)\};
$$
so
$R_{q,\sigma}(L)$ consists of $d$ rational functions of $\e$.
\begin{Proposition}
\label{eeppss}
With notations as above, if  $\re{q} > \k(L)$ and $C_{q,\sigma}(L)=\emptyset $
then no rational function from
$R_{q,\sigma}(L)$ has a pole at $\e =0$.
\end{Proposition}
\noindent {\bf Proof:} Each of the $d$-tuples
(\ref{dtup}) is a $d$-tuple of formal Taylor series (whose
non-constant terms are equal to 0). Any rational function from
$R_{q,\sigma}(L)$ can be represented by its formal power series.
If $C_{q,\sigma}(L)=\emptyset $, then every such series has to be
a Taylor series. \qed

\begin{rem}
\label{RemarkDispersion}
Let $q\in \sigma +\N$ as above.
Set $q' = \sigma + n + 1$ where
$n \in \Z$ is some integer for which
there exists no integer $i>n$
with $a_d(\sigma + i) = 0$. Suppose that $q-q' \geq 0$.
If $\Phi$ is a solution of $\hat{L}$, then
$\Phi(q),\ldots,\Phi(q+d-1)$ can be computed from
$\Phi(q'),\ldots,\Phi(q'+d-1)$
without dividing by $\e$.
This implies that  no rational function from
$R_{q',\sigma}(L)$ has a pole at $\e =0$.
\end{rem}

\section{Existence of a desingularization; algorithms}
\subsection{The sets
$C_{q,\sigma}(L)$ and $R_{q,\sigma}(L)$ in the case of apparent
singularities}
\label{lin}
We will use the following known result
(see  \cite[Proposition 4.4 and Theorem 4.5]{Im3}).

\begin{Theorem}
\label{Ram}
( Ramis \cite{R1}; Barkatou \cite{Bar}; Immink \cite{Im3})
The difference equation $L(y)=0$ admits
linearly independent
meromorphic solutions  $F _1,\dots ,F _d$ that are holomorphic in
the right half-plane $\re z> \k (L)$.
Moreover,
for some sufficiently large
integer $N$,
the sequences $\{F_j(n)\}_{n > N}$, $j=1,
\dots, d$, are linearly independent, i.e.,
\begin{equation} \label{**} \left\vert\begin{array}{ccc}

F_1(n)&\dots& F_1(n+d-1)\\
\vdots& & \vdots\\
F_d(n)&\dots& F_d(n+d-1)
  \end{array}\right\vert \neq 0
\end{equation}
for all $n >N$.
\end{Theorem}
As a consequence of Theorem \ref{Ram}
we get that if
$L$ is of the form (\ref{eq10}), then
for any complex number $q$ with $\re q$ large enough,
there exist
meromorphic solutions  $F _1,\dots ,F _d$ that are holomorphic in
the half plane $\re z > \k(L)$, such that
\begin{equation}
\label{***}
\left\vert\begin{array}{ccc}
F_1(q)&\dots& F_1(q+d-1)\\
\vdots& & \vdots\\
F_d(q)&\dots& F_d(q+d-1)
  \end{array}\right\vert \neq 0.
\end{equation}

\begin{Proposition}
\label{propn}
Let $\sigma $ be an apparent $t$-singularity of $L\in\C [z,E]$,
then $C_{q,\sigma}(L)=\emptyset $ for all $q\in \sigma +\N $ such
that  $\re q> \k (L)$.
\end{Proposition}
\noindent {\bf Proof:}
Suppose that
$q\in \sigma +\N $, $\re q> \k (L)$
and $C_{q,\sigma}(L)\neq \emptyset $.
Then for any $q'$ such that $q'-q\in \N$ we have
$C_{q',\sigma}(L)\neq \emptyset $, since if
$\Phi$ is a solution of $\hat{L}$, then
$\Phi(q'),\ldots,\Phi(q'+d-1)$ can be computed from
$\Phi(q),\ldots,\Phi(q+d-1)$
without dividing by $\e$. So we can assume that $\re q$
is larger than any preassigned real number.
By Definition \ref{defapp} the functions
$F_1,\dots , F_d$ mentioned in the consequence of Theorem \ref{Ram}
have no pole at $\sigma$.
If  $C_{q,\sigma}(L)\neq \emptyset$ then
by Proposition \ref{prop3}
there exists a non-trivial relation of the form
$u_0F_{q,0}+u_1F_{q+1,0}+\cdots +u_{d-1}F_{q+d-1,0}=0$ in $C_{q,\sigma}(L)$
and so the columns of (\ref{***})
are linearly dependent which is a contradiction, because inequality
(\ref{***}) has to be valid for all $q$ large enough.\qed


\begin{Proposition}
\label{propn''}
Let $\sigma $ be a $t$-singularity of $L\in\C [z,E]$,
then $C_{q,\sigma}(L)=\emptyset $ for all $q\in \sigma +\N $ such
that  $\re q> \k (L)$ if and only if $\sigma$ is apparent.
\end{Proposition}
\noindent {\bf Proof:}
It is obvious that if
$C_{q,\sigma}(L)=\emptyset $ then $\sigma $ is apparent. The converse
is true by Proposition \ref{propn}.\qed

As a consequence of
Propositions  \ref{eeppss}, \ref{propn}
we get
\begin{Proposition}
\label{propn'}
Let $\sigma $ be an apparent $t$-singularity of $L\in\C [z,E]$,
$q\in \sigma +\N $,  $\re q> \k (L)$.
Then no rational function from $R_{q,\sigma}(L)$
has a pole at $\e =0$.
\end{Proposition}
\subsection{$t$- and $l$- desingularizations}
\label{sectionds}
\noindent {\bf Algorithm {\tt t-desing}}. \\
{\bf Input:} $L \in \C[z,E]$ with non-zero $E^0$ coefficient. \\
{\bf Output:} A $t$-desingularization of $L$.
\begin{enumerate}
\item Let $a_0, a_d \in \C[z]$ be the trailing and leading
coefficient of $L$.
\item Let $n \in \N$ be the {\em dispersion} of $a_d$, $a_0$, which
is the largest integer
 such that $a_d$ has some root that equals $n$ plus
 some root of $a_0$.
If such $n \in \N$ does not exist then set $n$ equal 0.
\item Set
$L_2 := \frac{1}{a_0} L \in \C(z)[E]$.
\item For $i$ from 1 to $n$,
clear the $E^i$ coefficient of $L_2$ by subtracting $\frac{c}{a_0}
E^i L$ from $L_2$ where $c$ is the $E^i$ coefficient of $L_2$.
\item
Let $b_0 \in \C[z]$ be the least common multiple of the denominators
of the coefficients of $L_2$, and set $L_3 := b_0 L_2$.
\item
Compute $s,t \in \C[z]$ for which $s a_0 + t b_0 = {\rm gcd}(a_0,
b_0)$.  \item Output: $s L + t L_3$.  \end{enumerate}

Here is a Maple
implementation of {\tt t-desing}:
\begin{verbatim}
t_desing := proc(L, E, z)
   local a0, ad, n, L2, i, b0, L3;
   a0 := coeff(L,E,0); # trailing coefficient
   ad := lcoeff(L,E);  #  leading coefficient
   n := LREtools[dispersion](ad, a0, z, 'maximal');
   if not type(n,integer) then n := 0 end if;
   L2 := L/a0; # is in C(z)[E] with trailing coefficient 1
   for i from 1 to n do
      # Simplify L2 and clear its E^i coefficient:
      L2 := collect(L2, E, Normalizer);
      L2 := L2 - coeff(L2,E,i) * subs(z = z+i, L/a0) * E^i
   end do;
   # Multiply L2 by its denominator to obtain L3 in C[z,E]
   L3 := primpart(L2, E);
   b0 := coeff(L3, E, 0); # Equals the denominator of L2.
   gcdex(a0,b0,z,'s','t'); # gcd(a0, b0) = s*a0 + t*b0
   collect(s*L+t*L3, E, factor) # Output: s*L+t*L3 simplified
end proc:
\end{verbatim}
\begin{Theorem}
\label{tdes}
The algorithm {\tt t-desing} produces
a $t$-desingularization of $L$.
\end{Theorem}
\noindent {\bf Proof:}
Let $\sigma$ be an apparent $t$-singularity.
If $n$ is a non-negative integer
then there exists precisely one operator of the form
$$ L_2 = r_{n+d}(z) E^{n+d} + \cdots + r_{n+1}(z) E^{n+1} + 1  \in
\C(z)[E]$$
that is right-divisible by $L$  (recall that $d = \ord L$).
Algorithm {\tt t-desing} computes this operator $L_2$
where $n$ is the dispersion of $a_d$ and $a_0$.
If $u(z)$ is a solution of $L$, then it is also a solution of $L_2$
and so
$u(z) = -\sum_{i=1}^{d} r_{n+i}(z) u(z+n+i)$.
Plugging in $z=\e + \sigma$ we get
\[ u(\e+\sigma) = - \sum_{i=1}^{d}
r_{n+i}(\e+\sigma) u(\e+\sigma+n+i). \]
By taking the initial values $(1,0,0,\ldots,0)$,
$(0,1,0,\ldots,0)$, etc.,
for the $u(\e+\sigma+n+i)$ we obtain
the rational functions in the set $R_{\sigma+n+1,\sigma}(L)$
on the left-hand side of this equation, and the rational
functions $-r_{n+i}(\e+\sigma)$ on the right-hand side.
Since $\sigma$ is apparent, the rational functions in the set
$R_{\sigma+n+1,\sigma}(L)$ have no pole at $\e=0$ by
Proposition~\ref{propn'} when $n$ is sufficiently large (the
dispersion is large enough by Remark~\ref{RemarkDispersion}). Hence,
the $-r_{n+i}$ have no pole at $\sigma$.
Thus, $\sigma$ is not a root
of the denominators of $L_2 \in \C(z)[E]$. The least common multiple
of these denominators equals the trailing coefficient $b_0$ of $L_3$,
so $\sigma$ is not a root of $b_0$.
The trailing coefficient of the output is $s a_0 + t b_0 = {\rm gcd}(a_0,b_0)$.
Since this is a factor of $a_0$, the $t$-singularities of the output form
a subset of the $t$-singularities of $L$. And since this gcd divides
$b_0$, it follows that $\sigma$ is not a $t$-singularity of the output.
The same argument applies to every apparent $t$-singularity $\sigma$,
and hence the output is a $t$-desingularization of $L$.
\qed

Therefore the following theorem
(the main theorem for $t$-singularities)
is proven:
\begin{Theorem}
\label{mainth}
Every $L\inop$ is
$t$-desingularizable (in other words, there exists a
$t$-desingularization $\tilde{L}$ of $L$).
\end{Theorem}

\begin{Example}
For the operator 
\[
L\, = \, \left( 2\,z-1 \right)  \left( z-1 \right) {E}^{2}+ \left( 5\,z-1-9\,{z}^{2}+2\,{z}^{3} \right) E+z \left( 1+2\,z \right),
\]
algorithm {\tt t-desing} returns:
{
$1/3\, \left( -1+4\,z \right)  \left( 5\,z-1-9\,{z}^{2}+2\,{z}^{3} \right) {E}^{3}+ \left( {\frac {26}{3}}\,{z}^{2}-{\frac {43}{3}}\,z+11/3+{\frac {85}{3}}\,{z}^{3}\\
\mbox{}-18\,{z}^{4}+8/3\,{z}^{5} \right) {E}^{2} \\
\mbox{}+1/3\, \left( 7+4\,z \right)  \left( 5\,z-1-9\,{z}^{2}+2\,{z}^{3} \right) E+1$}
\end{Example}

The paper \cite{AvH99} gives a
proof that the so-called $\e$-criterion is
an alternative
a way to
decide if a desingularization exists.  The $\e$-criterion is based on
a construction similar to $R_{q,\sigma}(L)$. However,
in \cite{AvH99} a different definition of a desingularization was
used; the definition that we use in this paper is stricter (as
mentioned in Section \ref{intro}).

If the dispersion, the number $n$ in the algorithm,
is not positive then if follows from Theorem \ref{tdes} that no $t-$singularitiy is apparent. 
This fact already followed from the approach in \cite{AvH99},
see the summary of \cite{AvH99} given in Section~\ref{42}.

\begin{defi}
We call $\tilde{L}$ a {\em complete $t$-desingularization} of $L$
if $\tilde{L} \in \C [z,E]$
is right-divisible by
$L$ and its trailing coefficient is a non-zero constant.
If a complete $t$-desingularization of $L$ exists, then we say that
$L$ is {\em completely $t$-desingularizable}.
\end{defi}
\begin{Proposition}
\label{full}
A complete
$t$-desingularization of $L$ exists if and only if all
$t$-singularities of $L$ are apparent.
\end{Proposition}
\noindent {\bf Proof:}
Suppose a complete $t$-desingularization $\tilde{L}$ exists.
If $F(z)$ is a right-holomorphic solution of $L$,
then it is also a right-holomorphic solution of $\tilde{L}$.
Then $F(z)$ must be holomorphic since the trailing coefficient
of $\tilde{L}$ is a non-zero constant. Hence all $t$-singularities
of $L$ are apparent.
Conversely, if all $t$-singularities of $L$ are apparent then
a complete $t$-desingularization exists by Theorem \ref{tdes}. \qed

Algorithm {\tt t-desing} removes at least all apparent
$t$-singularities by Theorem \ref{tdes}. If a complete $t$-desingularization
exists, in other words, if all $t$-singularities can be removed,
then the above proposition shows that algorithm {\tt t-desing} will do so.
However, this does not imply that {\tt t-desing} always removes
as many $t$-singularities as possible:

\begin{Example}
\label{exar}
Let $L=(z+2)^2(z-1)^2E-(z+1)z(z-2)^2$.
The $t$-singularities of $L$ are $-1,0,2$, none of
which are apparent.
The application of {\tt t-desing} to $L$
gives $$-(z+3)(z+4)^2 E^3 + z(z-2)^2.$$
So one $t$-singularity disappeared even though no $t$-singularity
was apparent.
Note that the $t$-singularity $0$ can be removed as well:
the operator
$$
4(z+4)^2E^3-3z(z+3)(z+4)E^2+
3(z+2)(z-1)^2 E+ 2(z-2)^2
$$
is right-divisible by $L$. The $t$-singularity $z-2$ can not be
removed, because if all $t$-singularities could be removed then
all $t$-singularities would have to be apparent, and this is not
the case.
\end{Example}
We also implemented an algorithm that removes all singularities that
can be removed, by reducing this problem to a linear
algebra problem over the constants. This implementation is available
at: \\
\verb+http://www.math.fsu.edu/~hoeij/papers/desing/+ \\
and tends to produce nicer desingularizations than {\tt t-desing}.
We used {\tt t-desing} in this paper because it is shorter
than the linear algebra based desingularization algorithm, and
because the proof that all apparent
singularities can be removed is easier with {\tt t-desing}.

\begin{Example}
\label{ex2}
For the operator $L$
\begin{equation}
\label{fromex2}
(z-3)(z-2)E+z(z-1)
\end{equation}
we get $C_{4,1}(L)=C_{4,0}(L)=\emptyset$.
So $L$ must be completely $t$-desingularizable;
algorithm {\tt t-desing} returns:
${\frac {1}{72}}\, \left( 5\,z-6 \right)  \left( z-3 \right)  \left( z-2 \right) ^{2} \left( z-1 \right) {E}^{4}+$ \\
${\frac {1}{72}}\, \left( 108+106\,z+5\,{z}^{3}+39\,{z}^{2} \right)  \left( z-3 \right)
\mbox{} \left( z-2 \right) E \mbox{}+1$. \\
Note that nicer desingularizations are possible, for example $(E+1)^5$ of order 5,
or $-(z+1)E^4+(17z-29)E^3+(17z+56)E^2+(-z+5)E+1$ of order 4.
The operator $L$ from Example \ref{ex1}
is not completely desingularizable
because  $C_{4,-1}(L)\neq \emptyset$ there.
\end{Example}
For definiteness, we considered the trailing singularities.
If one already has an implementation of a $t$-desingularization algorithm,
then one can obtain an $l$-desingularization algorithm
by changing a small number of lines. However, one can avoid this duplication
of code because one can reduce $l$-desingularization to
$t$-desingularization
(and to get the algorithm {\tt l-desing})
with the following trick: interchange
the roles of the leading
and trailing coefficient by using the automorphism of $\C[z,E,E^{-1}]$
given by $z \mapsto -z$, $E \mapsto E^{-1}$.
This trick was used
in an implementation
\cite{Mit}
of
{\tt ds},
an old naive desingularization algorithm
\cite{AvH99}.
\subsection{An operator which is a desingularization of $L$
related to both
trailing and leading coefficients}
\label{both}
So far we considered mainly trailing apparent singularities.
Note that
Theorem \ref{mainth} is valid, mutatis mutandis, for leading apparent
singularities and a desingularization related to the leading
coefficient.
The prefix ``$lt$-'' indicates that we consider both leading
and trailing singularities, just like $l$- and $t$- indicate
leading and trailing.
The following theorem
generalizes Theorem \ref{mainth}.
\begin{Theorem}
\label{Th4}
(The main theorem.)
Any
$L\in \C[z,E]$ is $lt$-desingularizable.
\end{Theorem}
\noindent{\bf Proof :}
Let $R_t, R_l \in \C (z)[E]$, $L_t, L_l \in \C[z,E]$
be such that $L_t=R_tL$ (resp. $L_l=R_lL$)
is a $t$-(resp. $l$)-desingularization of $L$.
Consider the operator
\begin{equation}
\label{ml}
{\tilde L} = L_t + E^{m}L_l,
\end{equation}
where
$m = \max\{ 1, \ord L_t - \ord L_l+1\}$. It is clear that
${\tilde L}$ belongs to $\C[z,E]$, and ${\tilde L}$ is
right-divisible by $L$  since ${\tilde L}= (R_t+  E^{m} R_l)L$.
The operator ${\tilde L}$ is simultaneously a $t$- and
an $l$-desingularization of $L$ because its $t$-singularities are
$t$-singularities of $L_t$ and its $l$-singularities are $l$-singularities
of $L_l$.\qed

Analogously to the $t$-case
an operator $L$ is completely $lt$-desingularizable
if and only if any half-holomorphic
solution of $L(y)=0$ is holomorphic on $\C$.
In other words, $L$ is completely $lt$-desingularizable
if and only if all singularities (leading and trailing) are apparent.

We name {\tt desingboth} the algorithm given above.
Thus, the algorithm is given in the following manner:
\begin{itemize}
\item Use algorithms {\tt l-desing} and {\tt t-desing} for constructing
${L}_l$ and ${L}_t$.
\item Return ${\tilde L} = L_t + E^{m}L_l$, where $m = \max\{ 1, \ord L_t - \ord L_l+1\}$.
\end{itemize}

\begin{Example}
\label{ex5}
The operator
$(z-2)E-z$
has complete $lt-$desingularization $E^3-3E^2+3E-1$.
\end{Example}

\section{An application, previous work, and a conjecture}
\subsection{An application of desingularization}
\label{appds}
As we have mentioned, it is useful to have
a desingularization ${\tilde L}$
of $L$ for solving the continuation problem. Below we consider
another application of the desingularization.

Suppose that the sequence $u(0),u(1),\ldots$ satisfies the
following relation
\[
    (1+16n)^2 u(n+2) - (224+512n)u(n+1) - (n+1)(17+16n)^2 u(n) = 0
\]
which corresponds to $L = (1+16z)^2 E^2-(224+512z)E-(z+1)(17+16z)^2$.
Assume that $u(0),u(1) \in \Z$. By substituting $z=0,1,\ldots$ in the relation
we find:
\[
    u(2) = 289 u(0) + 224 u(1), \ \ \
    u(3) = 736 u(0) + 578 u(1), \ldots
\]
One sees that $u(2),u(3) \in \Z$ since we assumed that $u(0),u(1) \in \Z$.
The question is now the following: \\

\noindent
{\em Prove that $u(n) \in \Z$ for every nonnegative integer $n$.} \\

Each time we use $L$ to compute the next term $u(n+2)$ from the two previous
terms $u(n),u(n+1)$ we perform additions, multiplications,
and one division, namely by the {\em leading coefficient} of $L$ which is
$(1+16n)^2$. How to prove that this division does not cause $u(n+2)$ be
become a fraction?
This question becomes easy if we
find an $l$-desingularization w.r.t. the $l$-singularity
of $L$; if we use the algorithm {\tt l-desing}
(Sect. \ref{sectionds}),
then it
produces the following operator:
\[ L_l = E^3+(\frac{7}{2} z - \frac{81}{32})E^2 - (z+11)E
 - \frac{1}{32}(143+112z)(z+1). \]
Now $u(n+3)$ can be computed from $u(n),u(n+1),u(n+2)$ using $L_l$. This can
only introduce powers of 2 in the denominator of $u(n+3)$ because $L_l$
has only powers of 2 in the denominator and has leading coefficient 1.
So the denominators in the sequence must be powers of 2, but must
simultaneously be odd (and hence equal to 1) because the leading
coefficient $(1+16n)^2$ of $L$ is always odd.
Hence $u(n) \in \Z$ for all nonnegative integers $n$.

\subsection{Algorithm {\tt ds} in [1]}
\label{42}
%
%
In this section we will review algorithm {\tt ds} from \cite{AvH99} and
compare it with algorithm {\tt t-desing}.
Consider the operator
\begin{equation}
    \label{exL}
    L = (z-1)(z-2)(z+1)E^2+(z^5-3z^3+3z+2)E+z^2(z+2)
\end{equation}
The trailing coefficient has integer roots $z=0$ and $z=-2$.
We would like to decide if there is a $t$-desingularization of $L$,
where the $t$ refers to fact that we are only considering
the trailing coefficient $z^2(z+2)$ of $L$.
We will start with the largest integer root first,
so we first only consider the factor $z^2$.
The question is if there exists an operator
$R_t \in \C(z)[E]$ for which $R_t L$ is in $\C[z,E]$ and has
a trailing coefficient that is not divisible by $z^2$, or even
not divisible by $z$.
If such $R_t \in \C(z)[E]$ exists, then there must
exist an $R_t$ with the same property in $z^{-1} \C[z^{-1}][E]$. One
sees this by first replacing each coefficient of $R_t$ in $\C(z)$ by
its series expansion in $\C((z)) = \C[[z]][z^{-1}]$, and then throwing
away all non-negative powers of $z$ to end up with an element of $z^{-1} \C[z^{-1}]$.
So, if an operator $R_t$ with the desired property (that $R_t L$
is in $\C[z,E]$ and that the factor $z^2$ of the trailing coefficient
has disappeared) exists in $\C(z)[E]$, then
there exists an operator $R_t$ in $z^{-1} \C[z^{-1}][E]$ with the
same property. We can write this $R_t$ as
\begin{equation}
\label{eqRt}
    R_t = \sum_{i=0}^N r_i E^i, \ \ \ r_i = \sum_{j=-K}^{-1} c_{ij} z^j
\end{equation}
for some non-negative integers $N,K$ and some $c_{ij} \in \C$.

Suppose that there is a $c_{ij} \neq 0$ for some $j<-2$. Then take $i$
minimal with this property. Now $i$ can not be 0 because then the
trailing coefficient of $R_t L$ can not be in $\C[z]$.
The $E^i$ coefficient of $r_i E^i L$ is $r_i (z+i)^2(z+2+i)$ which
has a pole of order $>2$ at $z=0$ since $(z+i)^2(z+2+i)$ is not
divisible by $z$ (this is the reason we started with $z^2$ and not with $z+2$).
Then the $E^i$ coefficient of $R_t L$ has a pole of
order $>2$ as well (the pole in the $E^i$ term of $r_i E^i L$ can not cancel
against the $E^i$ terms of $r_0 E^0 L,\ldots,r_{i-1}E^{i-1}L$ since
those terms have pole orders $\leq 2$ by the minimality of $i$).
This means that $R_t L$ is not in $\C[z,E]$ which is a contradiction.
It follows that there can be no $j<-2$ for which $c_{i,j} \neq 0$
for some $i$, thus we can take $K=2$ without loss of generality.

Now take $N$ minimal with $r_N \neq 0$. Then the leading coefficient
of $R_t L$ is $r_N (z-1+N)(z-2+N)(z+1+N)$. This must be in $\C[z]$,
however, $r_N$ only has negative powers of $z$. It follows that
$(z-1+N)(z-2+N)(z+1+N)$ must be divisible by $z$. Hence $N$ can
be no greater than the largest integer root of the leading coefficient,
which is 2. So we can take $N=2$.
In the more general situation where we want to eliminate
several roots of the trailing coefficient at the same time one can
use the same argument to show that one may assume without loss of generality
that the order of $R_t$ is bounded by the {\em dispersion}
(the largest root difference in $\Z$) of the leading
and trailing coefficient of $L$.

We now see that if $R_t$ exists, then we may assume it to be of the
form in equation~(\ref{eqRt}), and from the preceding we see that we
may also assume $N=2$ and $K=2$.
This turns the problem
into a finite dimensional system of linear equations for the $c_{ij}$.
Solving this system decides whether or not the factor $z^2$ can be removed,
and if so, how to do this. We find the following solution
\[ R_t = \frac{-1}{12z}E^2+(\frac{5}{9z} - \frac2{3z^2})E + \frac1{z^2}. \]
Then $R_t L$ is in $\C[z,E]$ and has trailing coefficient $z+2$.
Algorithm {\tt ds} in \cite{AvH99} finds this $R_t$ in a slightly
different way. Write $R_t = \sum r_i E^i = z^{-2} \sum \tilde{r}_i E^i$
where $\tilde{r}_i = z^2 r_i \in \C[z]$. We can now work modulo $z^2$,
so we can view $\tilde{r}_i$ as an element of $\C[z]/(z^2)$.
Now equate $(\sum \tilde{r}_i E^i)L$ to zero modulo $z^2$.
This leads to a system of linear equations for the $\tilde{r}_i$. This
system is already in a triangular form so it can be solved
quickly by applying Gaussian elimination over $\C[z]/(z^2)$.
If we set $\tilde{r}_0 := 1$ then one first finds
$\tilde{r}_1 = (5/9)z - 2/3$. This involves a division in $\C[z]/(z^2)$,
namely one divides by $a_0(z+1)$
where $a_0(z) = z^2(z+2)$ is the trailing coefficient of $L$.
After that one computes $\tilde{r}_2 =  -z/12$, which involves a division
by $a_0(z+2)$.

After eliminating $z^2$ we can apply a similar process to the factor $z+2$ in
the trailing coefficient. Then one obtains linear equations over $\C[z]/(z+2)$
instead of over $\C[z]/(z^2)$ and one can proceed along the same lines,
see \cite{AvH99} for details. It
turns out that $z+2$ can be removed as well.

Algorithm {\tt ds} tries to remove one root from the trailing
coefficient at a time.
If a root can not be removed, that is, if we encounter a non-apparent
$t$-singularity, then algorithm {\tt ds} stops, and later roots
will not be removed even if some of them correspond to apparent
$t$-singularities.
The algorithm {\tt t-desing} presented in this paper has the
advantage that it removes all apparent $t$-singularities,
even if there are non-apparent $t$-singularities between them.
Furthermore, it is shorter than algorithm {\tt ds} and does not need to compute with
roots of $a_0$.
Algorithm {\tt  t-desing} is also simpler than differential desingularization (see Appendix) since
it does not need to know which singularities are apparent.

\subsection{A conjecture}

Let $L \in \Q[z,E]$ with leading and trailing coefficient $a_d(z)$
and $a_0(d)$. Let $N_1$ be an integer larger than any integer root
of $a_d(z-d)a_0(z)$ and let $N_2$ be an integer smaller than any
integer root of $a_d(z-d)a_0(z)$. Then there exist $d$ linearly
independent sequences $u_1,\ldots,u_d: N_1 + \N \rightarrow \Q$
and $d$ linearly independent sequences $v_1,\ldots,v_d: N_2 - \N \rightarrow \Q$
that satisfy the recurrence given by $L$.
Let $P$ denote the set of all prime numbers $p$ that occur as a factor
in a denominator of a number that appears in a sequence $u_i$ or $v_i$.
So $p$ is not in $P$ if and only if every entry of every $u_i$ and $v_i$
does not have $p$ in the denominator.
We call $L$ {\em smooth} if the set $P$ is finite. Note that the
set $P$ may depend on the choices of the $u_i$ and $v_i$, but whether $P$
is finite or not only depends on $L$.
\\


\noindent {\bf Conjecture:} $L$ is smooth if and only if
$L$ is completely $lt-$desingularizable. \\

The conjecture relates analytic properties to number theoretic properties,
namely it states that sequence solutions (where we consider
sequences that extend to the right as well as sequences that extend to the
left) have only finitely many primes in the denominators iff all
half-holomorphic solutions are holomorphic.

If 
$L$ has a complete $lt$-desingularization $\tilde{L}$ then we can extend
sequences $u_i$ to the right and sequences $v_i$ to the left with
$\tilde{L}$. This will only introduce
finitely many primes in denominators (so $L$ is smooth)
since $\tilde{L}$ has constant
leading and trailing coefficient. This
shows that the conjecture is true in one direction.

\section *{Appendix: Differential Case}
\label{dffrntl}
In this appendix we prove that the statement in Theorem \ref{mainth}
is also true in the differential case. The proof is essentially a
desingularization algorithm (we name it $\partial$-{\tt desing}).
An implementation can be found at:  \\
\verb+http://www.math.fsu.edu/~hoeij/papers/desing/+ \\
However, the results in this appendix are not new; more general results were
given by Tsai in \cite{Tsai}
(a desingularization
is one of the elements in the output of Tsai's Weyl Closure algorithm).
We include this appendix for completeness and because the proof is
short since
the result is less general than \cite{Tsai}.

Let $\partial$ the differentiation $d/dz$.
If $L\in\C[z,\partial]$ and the coefficients of $L$
do not have a non-constant common factor,
then a {\em singular point} or {\em singularity} of $L$
is a zero of
the leading coefficient of $L$. If $L\in \C(z)[\partial]$
and $L$ is monic
\begin{equation}
\label{ap1}
L=\partial ^n+a_{n-1}(z)\partial ^{n-1}+\cdots +a_0(z),
\end{equation}
then singularities are  poles of  $a_i$'s. We will consider
differential operators in the form (\ref{ap1}).
Definitions of {\em regular
singularity} and, resp., {\em irregular singularity} of a given
operator
$L$ of the form (\ref{ap1})
can be found, e.g., in \cite{Ince}. A point that is not
a singularity (neither regular nor irregular) is {\em ordinary}.

\begin{defi}
Let $L\in \C(z)[\partial]$.
A singularity $p$ of $L$
is {\em apparent} if
there exists an open set $U$ with $p \in U$ and a basis of solutions
of $L(y)=0$ that are holomorphic on $U$.
\end{defi}


\begin{defi}
\label{mvhd3}
Let $L \in \C (z)[\partial]$ be monic.
Then $\tilde{L}$ is called a {\em desingularization} of $L$ if it has $L$ as a
right-hand factor, and if for every $p \in \C$, if $p$ is either
an ordinary point of $L$ or
an apparent singularity of $L$, then $p$ is an ordinary point of
$\tilde{L}$.
\end{defi}
\begin{Proposition}
\label{mvhl1}
Let $L \in \C (z)[\partial]$ be monic, have order $n$, and let $p \in
\C$.  The following statements are equivalent.
\begin{enumerate}
\item[a)]  There exist an open set $U$ with $p \in U$ and a basis of
solutions $y_1,\ldots,y_n$ of $L$, holomorphic on $U$, for which
$y_i$ vanishes at $p$ with order $i-1$.

\item[b)]
$p$ is not an irregular singularity of $L$,
 the local exponents at $p$ are $0,1,\ldots,n-1$, and the
    formal solutions of $L$ at $p$ are in $\C [[z-p]]$.

\item[c)]  $p$ is not a pole of any of the coefficients in $\C (z)$
of $L$ (i.e., $p$ is an ordinary point of $L$).
\end{enumerate}
\end{Proposition}
\noindent {\bf Proof:} a)
$\Rightarrow$ b) is clear; c) $\Rightarrow$ a) is Cauchy's theorem; a
proof for b) $\Rightarrow$ c) can be found in \cite[Lemma
9.2]{vh97}. \qed
\begin{Proposition}
\label{mvhl2}
Let $L \in \C (z)[\partial]$
be monic, have order $n$, and let $p \in \C$.
The following statements are equivalent.
\begin{enumerate}
\item[a)]
$p$ is either an ordinary point or an apparent singularity of $L$.
\item[b)]  $L$ has $n$ linearly independent solutions in $\C [[z-p]]$
    (this is equivalent
    to: $p$ is not an irregular singular point of
$L$, the local exponents are
    non-negative integers and the formal solutions at $p$ do not
    involve logarithms, see \cite{vh97} for more details).

\item[c)]  There exists a monic operator $\tilde{L} \in
\C (z)[\partial]$ that has $L$
    as a right-hand factor such that $p$ is an ordinary point of
$\tilde{L}$.  \end{enumerate} \end{Proposition}
\noindent {\bf Proof:}
a) $\Rightarrow$ b) is clear; c) $\Rightarrow$ a) follows
from Cauchy's theorem.  For b) $\Rightarrow$ c), let $m$ be the
highest local exponent.  Consider the set $V$ of all integers $0 \leq
i \leq m$ for which $i$ is not an exponent of $L$ at $p$.  Let $L_1$
be an operator with the following as basis of solutions:  $L((z-p)^i
)$,  $i \in V$.  Let $\tilde{L} = L_1 L$.  Then $\tilde{L}$ satisfies
part a) of Proposition \ref{mvhl1}. \qed

\begin{rem}
In general, solutions in $\C [[z-p]]$
need not be convergent, but b) $\Rightarrow$ a)
of
Proposition
\ref{mvhl2}
says that if {\em all} formal solutions of $L$ are in
$\C [[z-p]]$
then they are automatically convergent.
\end{rem}

%
%

\begin{Theorem}
\label{mvht1}
A desingularization always exists.
\end{Theorem}

\noindent {\bf Proof:}  Let $A$ be the set of all apparent
singularities $p \in \C$. For $p \in A$, let $m(p)$ be the
highest exponent at $p$. Let $m$ be the maximum of all $m(p)$,
and let $e(p) \subset \{0,1,\ldots,m\}$  be the set of
exponents of $L$ at $p$.

Take $p \in A$ for which $m = m(p)$.
Suppose a desingularization $\tilde{L}$ exists. The exponents of $L$ at $p$ must
be a subset of the exponents of $\tilde{L}$ at $p$ because $L$ is a right-hand factor
of $\tilde{L}$.  Hence $m$ is an exponent of $\tilde{L}$ at $p$, but since $p$ is
a ordinary point
of $\tilde{L}$ it follows that $0,1,\ldots,m$ are exponents of $\tilde{L}$
at $p$ as well, so the order of $\tilde{L}$ must be at least $m+1$.

We will now show that a desingularization $\tilde{L}$ of order
$m+1$ exists. First, construct polynomials
$y_1, \ldots, y_{m+1-n}
\in \C [z]$ such that for every $p \in A$ and every $i \in
\{0,1,\ldots,m\} \setminus e(p)$ there is precisely one $y_j$
that vanishes at $p$ with order $i$. Let $L_1$ be the monic
operator whose solutions are spanned by
$L(y_1),\ldots,L(y_{m+1-n})$. Then $L_1 L$ satisfies condition a)
of Proposition~\ref{mvhl1} at every $p \in A$. Thus, $p$ is
an ordinary point of
$L_1 L$ for every $p \in A$.  However,  $L_1 L$ need not satisfy
the definition of a desingularization because we may have created
new apparent singularities. Now $L_1 L$ has order $m+1$ and so
$L_1$ has order $m+1-n$. Write $L_1 = \partial^{m+1-n} + \sum a_i
\partial^i$ with $a_i \in \C (z)$. Take $b_i \in \C (z)$ as
follows:  $b_i$ has no poles in $\C \setminus A$, and $b_i - a_i$
must vanish at every $p \in A$ with order at least $M(p)$ for
some integers $M(p)$.  We take these integers $M(p)$ high enough
to make sure that $L_2 L$ is an element of $\C (z)[\partial]$
whose coefficients in $\C (z)$ have no pole at any $p \in A$,
where $L_2 = \sum (b_i - a_i)\partial^i$.  Now let $L_3 =
\partial^{m+1-n} + \sum b_i \partial^i$ and $\tilde{L} = L_3 L$.
Since the $b_i$ have no poles in $\C \setminus A$, we see that
every
ordinary point in $\C$ of $L$ is an ordinary point of $L_3$
and hence an ordinary point of $\tilde{L}$.
The operators $L_1 L$
and $L_2 L$ have coefficients that do not have poles at any $p \in
A$, hence the same is true for $\tilde{L} = L_1 L + L_2 L$, and
since $\tilde{L}$ is monic it follows that every
$p \in A$ is an ordinary point of
$\tilde{L}$.\qed
\begin{Example}
\label{dfrncl}
Let $L$ be the monic operator with $z\,{\rm cos}(z)$ and $z\,{\rm sin}(z)$
as basis of solutions.
Then $L$ has one singularity in $\C$, namely at $z=0$, which is
an apparent singularity with local exponents $1$ and $2$.
Thus, to desingularize $L$ we must add a solution with exponent 0.
Take $y_1 = z^0 = 1$ and compute $L(y_1)$. The explicit form of $L$ is
\[ L = \partial^2 - \frac2z \partial + 1 + \frac2{z^2} \]
so we find $L(y_1) = 1+2/z^2$. Now let $L_1$ be the monic operator with $1+2/z^2$
as a basis of solutions, then
\[ L_1 = \partial + \frac4{z(z^2+2)}. \]
Multiplying $L$ on the left by $L_1$ adds a solution (namely $y_1=z^0$) to $L$
with the missing exponent 0.
Then $L_1 L$ satisfies part b) of Proposition~\ref{mvhl1}. Hence $z=0$ is
a regular point of $L_1 L$.
Unfortunately, $L_1$ introduces new singularities, namely at $z^2+2=0$.
We will illustrate how to remedy this with a
truncated Laurent series at $z=0$.
We have \[
L_1 = \partial + \frac4{z(z^2+2)}
= \partial + 2z^{-1} - z + \frac12 z^3 - \frac14 z^5 + \cdots \]
Since the highest power of $z$ in the denominator of $L$ is $z^2$, we see that
if we change any $z^i$-term in $L_1$ with $i \geq 2$ then this
can not introduce poles at $z=0$ in the coefficients of $L_1 L$.
If we remove these terms from $L_1$ we get an operator
$L_3 := \partial + 2z^{-1} - z$ for which the coefficients
of $L_3 L$ still have no pole at $z=0$. Since $L$ and $L_3$ have no other
singularities in $\C$ apart from $z=0$, we see that $\tilde{L} := L_3 L$ has
no other singularities in $\C$ either, so $\tilde{L}$ must be in $\C[z,\partial]$
and indeed:
\[ \tilde{L} = \partial^3 - z\partial^2 + 3\partial - z. \]
\end{Example}

Note that this desingularization process generally makes the singularity at $z=\infty$ worse
in the sense that the highest slope in the Newton polygon at $z=\infty$ is higher for $L_3$
(and hence for $\tilde{L}$) than it
is for $L$. In
Example \ref{dfrncl}
this is unavoidable if we want
$\tilde{L}$ to have minimal order, although for this particular example
a desingularization $\tilde{L} = (\partial^2+1)^2$
of non-minimal order exists whose singularity at $z=\infty$
is not worse than that of $L$.
Our implementation always computes a desingularization of minimal order
but it is not optimal in the sense that it makes no effort to avoid
making the singularity at $z=\infty$ worse even in examples
where this could be done.

\begin{defi}
We call $\tilde{L}$ a {\em complete desingularization}
of $L$
if $\tilde{L} \in \C [z,\partial]$ has $L$ as
right-hand factor and is monic
with respect to $\partial$.
\end{defi}

Theorem \ref{mvht1}
implies that a complete desingularization
of $L$ exists if and only if all singularities of $L$
in $\C$ are apparent singularities.

\end{document}